\documentclass[12pt]{article}

\usepackage{hyperref}
\font\smallit=cmti10
\font\smalltt=cmtt10
\font\smallrm=cmr9

\begin{document}

\begin{center}
{\bf TAMING THE WILD IN IMPARTIAL COMBINATORIAL GAMES}
\vskip 20pt
{\bf Thane E. Plambeck }\\
{\smallit 2341 Tasso Street, Palo Alto, California 94301}\\
{\tt thane@best.com}\\ 
\href{http://www.plambeck.org/oldhtml/mathematics/games/misere}{\tt http://www.plambeck.org/oldhtml/mathematics/games/misere}\\
\end{center}
\centerline{\smallit Received: 4 Dec 2004, Revised: 8 Jan 2005, Accepted: 29 May 2005, Published: }
\vskip 30pt

\centerline{\bf Abstract}

\noindent
We introduce a {\em misere quotient semigroup} construction 
in impartial combinatorial game theory, and argue
that it is the long-sought natural generalization of the normal-play 
Sprague-Grundy theory to misere play.  Along the way, we illustrate how to
use the theory to describe complete analyses
of two wild taking and breaking games.   

\pagestyle{myheadings}
\markright{\smalltt INTEGERS: \smallrm ELECTRONIC JOURNAL OF COMBINATORIAL NUMBER THEORY \smalltt x (200x), \#Axx\hfill}

\thispagestyle{empty}
\baselineskip=15pt
\vskip 30pt 

\section{\normalsize Introduction}

On page 146 of {\em On Numbers and Games}, in Chapter 12, ``How to Lose When You Must,''
John Conway writes:

\begin{quotation}
{\em Note that in a sense, [misere] restive games are {\em ambivalent Nim-heaps}, which
choose their size ($g_0$ or $g$) according to their company.  There are
many other games which exhibit behaviour of this type, and it would
be very interesting to have some general theory for them.}
\end{quotation}

This paper provides such a general theory, cast in the language of 
commutative semigroups.  We have two goals:
\begin{itemize}
\item Generalize the normal-play Sprague-Grundy theory of impartial games to misere play. 
\item Describe complete winning strategies for particular wild misere games.
\end{itemize}

We introduce a {\em quotient semigroup} 
structure on the set of all positions of an impartial game with fixed rules.
The basic construction is the same for both normal and misere play.  In normal play, 
it leads to the familiar Sprague-Grundy theory.  In misere play, when applied to 
the set of all sums of positions played according to 
a particular game's rules, it leads to a quotient of a free commutative semigroup by
the game's {\em indistinguishability congruence}. 
Playing a role similar to the one that {\em nim sequences}
do for normal play, mappings from single-heap
positions into a game's misere quotient semigroup 
succinctly and necessarily 
encode all relevant information about its best misere play.  Studying examples in 
detail, we'll see how wild misere games that 
involve an infinity of ever-more complicated canonical forms amongst their position sums may nevertheless 
possess a relatively simple, even finite misere quotient.  
Many previously unsolved wild misere games have now been resolved \cite{mg} using such techniques.

\section{\normalsize Prerequisites}
We're going to assume that our reader is already familiar with the theory of 
normal- and misere-play
impartial combinatorial games as presented in \cite{ww1}, \cite{ww2}, or \cite{onag}.
For basic concepts and results cited from commutative semigroup theory, we 
follow \cite{cp}.

Specifically, our reader should be familiar with the following: 
the abstract definition of an {\em impartial game}; the
convention that a {\em sum} of games is played {\em disjunctively}; the
difference between the {\em normal play} (last player winning) and {\em misere play}
(last player losing)
{\em play conventions};
the rules of the game of {\em Nim}; 
the {\em Sprague-Grundy theory}, including {\em nim equivalents} (ie, nim-heaps $\ast k$), 
{\em nim-addition}, and the {\em mex rule}; the idea that each impartial game
has a deterministic {\em outcome class} that describes it as either an {\em P-position}
 (previous player winning in best play) 
or an {\em N-position} (next player winning); the notion of
{\em canonical forms} for normal and misere play games, and how to compute them; the
{\em game code} notation for specifying the rules of a {\em taking and breaking game}, 
and related concepts.   

For misere play, we're going to additionally assume that the reader knows 
what {\em genus symbols} are, how to compute them, their relation to correct play of 
{\em misere Nim},
and the role they play in the {\em tame-wild} distinction.  We'll try to limit our 
dependence on these latter concepts, however.


\vskip 30pt

\section{\normalsize Misere taking and breaking games as semigroup quotients}

Suppose $\Gamma$ is a taking and breaking game whose rules have been fixed in advance.  
The reader is invited to think of $\Gamma$ as standing for Nim, or Dawson's Chess,
or Kayles, or any other impartial game that can be played using heaps of beans.

\subsection{\normalsize Heap alphabets}
Let $h_i$ be a distinct, purely formal symbol for each \mbox{$i \geq 1$}.
We will call the set $H = \{h_1,h_2,h_3,\ldots \}$ the {\em heap alphabet}.
A particular symbol $h_i$ will sometimes be called a {\em heap of size $i$}.

The notation $H_n$ stands for the subset $H_n = \{h_1, \ldots, h_n\} \subseteq H$ 
for each $n \geq 1$.

\subsection{\normalsize Game positions} 
Let $\mathcal{F}_H$ be the free commutative semigroup 
on the heap alphabet $H$.  The semigroups
$\mathcal{F}_H$ and ${\mathcal{F}}_{H_n}$ include an identity
$\Lambda$, which is just the empty word.

There's a natural correspondence between 
the elements of $\mathcal{F}_H$ and the set of
all position sums  
of a taking and breaking game $\Gamma$.  In this correspondence, a finite sum of 
heaps of various sizes is written multiplicatively using corresponding elements of the
heap alphabet $H$.
For example, a position with two heaps of
size six, and one each of sizes three and two would correspond to the product
\[h_2h_3{h_6}^2.\]

This multiplicative notation for sums makes it convenient to take the convention that the
empty position $\Lambda=1$.  It corresponds to the {\em endgame}---a position with no options.

\vspace{0.1in}
The following definition and simple lemma are of the utmost importance to us.

\subsection{\normalsize Indistinguishability} 
\label{indistinguishability}
Fix the rules and associated {\em play convention} (normal or misere) 
of a particular taking and breaking game $\Gamma$.
Let $u, v \in \mathcal{F}_H$ be game positions in $\Gamma$.   We'll say that
{\em $u$ is indistinguishable\footnote{The indistinguishability definition originates in Conway's (\cite{onag} pg 147ff) discussion
of misere canonical forms, where 
the elements $w$ were taken instead from the set of {\bf all} impartial misere combinatorial
games.  We're interested in indistinguishability only over
the smaller set (and subsets of) $\mathcal{F}_H$, the set of all sums that actually arise 
in the game $\Gamma$ we're studying.  
Although the difference may seem slight, it is in fact crucial to
the success of the methods described in this paper.  See \S \ref{canonical} for more discussion.}
from $v$ over $\mathcal{F}_H$},
and write the relation $u \ \rho \ v$, if for every element $w \in \mathcal{F}_H$,
$uw$ and $vw$ are either both $P$-positions, or are 
both $N$-positions.

{\bf Lemma 1}
The relation $\rho$ is a congruence on $\mathcal{F}_H$.

{\it Proof.}
We must show that $\rho$ is a reflexive, symmetric, transitive, and compatible 
relation on $\mathcal{F}_H$.  It's easy to see that the indistinguishability
definition ensures that $u \ \rho\ u$ is always true, and
that $u \ \rho \ v$ implies $v \ \rho \ u$.  So $\rho$ is reflexive and symmetric.
\begin{itemize}
\item{$\rho$ {\em is transitive}}: Suppose $u \ \rho \ v$ and $v \ \rho \ s$.  Since $u \ \rho \ v$, for every
choice of $w \in \mathcal{F}_H$, $uw$ and $vw$ have the same outcome.  
Since $v \ \rho \ s$, $vw$ and $sw$ have the same outcome.  So $uw$ and $sw$ have the same outcome,
ie, $u \ \rho \ s$.
\item{$\rho$ {\em is compatible}}: We need to show that  $u \ \rho \ v$ implies 
$uw \ \rho \ vw$ for every $w \in \mathcal{F}_H$.  So suppose that $s$ is an arbitrary element of $\mathcal{F}_H$.
We need to show that $uws$ and $vws$ have the same outcome.  But if we let $w' = ws$, we can use that
fact that $u \ \rho \ v$ to conclude that $uw'$ and $vw'$ have the same outcome.  So 
$uw \ \rho \ vw$.
\end{itemize}

We come now to the central object of our study.

\subsection{\normalsize The quotient semigroup}
\label{qs}
Suppose the rules and play convention 
of a taking and breaking game $\Gamma$ are fixed, and let $\rho$ be the indistinguishability
congruence on $\mathcal{F}_H$, the free commutative semigroup of all positions in $\Gamma$.
The {\em indistinguishability quotient} $\mathcal{Q} = {\mathcal Q}(\Gamma)$ is the commutative semigroup
\[Q = \mathcal{F}_H/\rho.\]

Notice that the indistinguishability quotient can be taken with respect to either play convention (normal
or misere).  The details of the indistinguishability congruence then determine the structure of the
indistinguishability quotient.  Since the word ``indistinguishability'' is quite a mouthful, 
we prefer to call $\mathcal{Q}$ the {\em quotient semigroup} of $\Gamma$.

When $\Gamma$ is a normal play game, 
its quotient semigroup $\mathcal{Q}= {\mathcal Q}(\Gamma)$ is more than just a semigroup.  The Sprague-Grundy theory says that 
it is always a {\em group}.   It's isomorphic
to a direct product of a (possibly infinite) set of ${Z}_2$'s (cyclic groups of order two).  
If $u$ is a position in $\mathcal{F}_H$ with normal play nim-heap equivalent $\ast k$,
the members of a particular congruence class $u \rho \in \mathcal{F}_H/\rho$ 
will be precisely all positions that have normal-play nim-heap equivalent $\ast k$.
The identity of $\mathcal{Q}$ is the congruence
class of all positions with nim-heap equivalent $\ast 0$.
The ``group multiplication'' corresponds to nim addition.  
We won't have much more to say about such normal play quotients in this paper.
Instead, we'll be almost exclusively interested in {\em misere quotient semigroups}. 

For misere play, the quotient structure is a {\em semigroup}.  Surprisingly, it's often a finite object,
even for a game that has an infinite number of different canonical forms occurring amongst its sums.  

Do the elements of a particular congruence class all have the same outcome?  Yes.
Each class can be thought of as carrying a big stamp labelled ``P'' (previous player
wins in best play for all positions in this class) or ``N'' (next player wins).  
In normal play, there's only one equivalence class labelled ``P''---these are the positions with
nim heap equivalent $\ast 0$.   In misere play, for all but the
trivial games with one position $\ast 0$, or two positions $\{\ast 0, \ast 1\}$, there is always more than 
one ``P'' class---one corresponding to the position $\ast 1$, and at least one more, corresponding to the position $\ast 2 + \ast 2$.

It's time to look at a concrete example of the quotient semigroup of a wild misere game.

\section{\normalsize  How to lose at {\bf 0.123}}
\label{howtolose}
The octal game {\bf 0.123} can be played with counters arranged in heaps.
Two players take turns removing one, two or three counters from a heap,
subject to the following conditions:
\begin{enumerate}
\item Three counters may be removed from any heap;
\item Two counters may be removed from a heap, but only if it has more than two counters; and
\item One counter may be removed only if it is the only counter in that heap.
\end{enumerate}

In {\it normal play} of {\bf 0.123}, the 
last player able to make a legal move is declared the winner.  
In normal play, each heap size reduces to a nim-heap.
The resulting nim sequence\footnote{See {\it Winning Ways} \cite{ww2}, Vol I, Chapter 4, pg 87, ``Other Take-Away Games;'' also Table 7(b),
pg 104.} is periodic of length 5, beginning at heap 5.  See Figure \ref{0123normal}.

\begin{figure}[h]
\begin{center}
\begin{tabular}{c|ccccc}
+   & 1 & 2 & 3 & 4 & 5 \\ \hline
0+ & 1 & 0 & 2 & 2 & 1 \\
5+ & 0 & 0 & 2 & 1 & 1 \\
10+ & 0 & 0 & 2 & 1 & 1 \\
15+ & $\cdots$ & & & &
\end{tabular}
\caption{\label{0123normal} Normal play nim heap equivalents for ${\bf 0.123}$.}
\end{center}
\end{figure}

In {\it misere play}, the last player to make a legal move is 
declared to be the {\it loser} of the game.  

Taking our notation from {\it Winning Ways} (\cite{ww2}, Vol II, Chapter 13, ``Survival in the Lost World''), 
we exhibit the {\em genus sequence} of misere {\bf 0.123} in Figure \ref{0123misere}.  This sequence is
also periodic of length 5.  See Figure \ref{0123misere}.

\begin{figure}[h]
\begin{center}
\begin{tabular}{c|ccccc}
+   & 1 & 2 & 3 & 4 & 5 \\ \hline
0+ & 1 & 0 & 2 & 2 & 1 \\
5+ & $0^{02}$ & 0 & $2^{1420}$ & $1^{20}$ & 1 \\
10+ & $0^{02}$ & 0 & $2^{1420}$ & $1^{20}$ & 1 \\
15+ & $\cdots$ & & & &
\end{tabular}
\caption{\label{0123misere}{\sf G}*-values of {\bf .123}}
\end{center}
\end{figure}In Figure \ref{0123misere}, an entry that is a simple integer (0, 1, or 2) represents
that the game at that position has a misere canonical form identical to
a misere nim heap of the corresponding size.  
The 
genus symbols\footnote{See \cite{ww2}, Vol II, pg 422, ``Animals and Their Genus.'' In \cite{ww1}, see Vol I, pg 402.} 
of the nim heaps 
that occur in Table 2 are
\begin{eqnarray*}
0 & = & 0^{1202020\cdots} = 0^{120} \\
1 & = & 1^{0313131\cdots} = 1^{031} \\
2 & = & 2^{2020202\cdots} = 2^{20}.
\end{eqnarray*}
In misere play of {\bf 0.123}, the first non-nim-heap
occurs at the six-counter heap.  It 
is the game $h_6=2_{+} = \{2\}$.   The eight-counter heap
is $h_8 = \{2_{+}, 1\}$, and the nine-counter heap is $h_9 = \{h_8, 0\}$.
Unlike $h_6$,
the latter two positions are {\em wild}---their genera match the genus of
no misere Nim position.
Although the subsequent canonical forms of the games occurring at
heap sizes = 1, 3, and 4 (modulo 5) are not identical to
$h_6$, $h_8$, and $h_9$, respectively, their respective genera do repeat,
as indicated in Figure \ref{0123misere}.

The information in Figure \ref{0123misere} is sufficient to determine outcome classes for
{\em single heap} misere {\bf 0.123} positions, and also sums of a single heap with
arbitrary numbers of nim heaps of size one or two (via the genus symbol exponents).  
To capture information
about the best play of {\em all} 
misere {\bf 0.123} positions (ie, an arbitrary sum of arbitrarily-sized heaps), we change our viewpoint entirely---we write down a semigroup
presentation for its misere quotient $\mathcal{Q}$.

\subsection{\normalsize  The misere quotient semigroup $\mathcal{Q}_{\bf 0.123}$}

\begin{figure}[h]
\begin{center}
\begin{tabular}{c|ccccc}
  & 1 & 2 & 3 & 4 & 5 \\ \hline 
 0+  & $x$  & $e$  & $z$  & $z$  & $x$  \\  
 5+  & $b^2$  & $e$  & $a$  & $b$  & $x$  \\  
 10+  & $b^2$  & $e$  & $a$  & $b$  & $x$  \\  
 15+  & $\cdots$  &   &   &   &   \\  
\end{tabular}
\end{center}
\caption{\label{0123-semi} Semigroup identifications for single heaps in misere {\bf 0.123}}
\end{figure}

{\bf Theorem 1.}
\label{0123}
The misere quotient $\mathcal{Q}(\Gamma)$ of the wild octal game $\Gamma = {\bf 0.123}$ is isomorphic to a
twenty-element commutative semigroup $\mathcal{Q}_{\bf 0.123}$ with identity $e$ that is 
presented by the following generators and relations:
\[\{x,z,a,b \ | \ x^2=a^2=e, \ z^4=z^2, \ b^4=b^2, \ abz=b, \ b^3x=b^2, \ z^3a=z^2 \}.\]

We won't have the tools in place to {\em prove} that $\mathcal{Q}$ and $\mathcal{Q}_{\bf 0.123}$ 
are isomorphic for several more sections 
(see section \ref{proving} if you can't wait).   Our immediate goal is to
take a closer look at $\mathcal{Q}_{\bf 0.123}$ itself.  

The twenty elements of $\mathcal{Q}_{{\bf 0.123}}$ are
\[\{e, x, z, a, b, xz, xa, xb, z^2, za, zb, b^2, xz^2, xza, xzb, xb^2, z^3, zb^2, xz^3, xzb^2\},\]
and they can be partitioned into fifteen pairwise distinguishable N-position elements 
\[\{e, z, a, b, xz, xb, za, xz^2, xza, xzb, xb^2, z^3, zb^2, xz^3, xzb^2\},\]
and five pairwise distinguishable P-position elements
\[\{x, xa, b^2, z^2, zb\}.\] See Figure \ref{0123action}.

\begin{figure}[h]
\begin{center}
\begin{tabular}{c|c|c|c||cccc}
\# & element & genus & outcome & x  & z  & a & b \\ \hline
1  & $e$     & $0^{120}$  & N & 2  & 3  & 4  & 5 \\
2  & $x$     & $1^{031}$  & P & 1  & 6  & 7 & 8 \\
3  & $z$     & $2^{20}$   & N & 6  & 9  & 10  & 11 \\
4  & $a$     & $2^{1420}$ & N & 7  & 10  & 1  & 11  \\
5  & $b$     & $1^{20}$   & N & 8  & 11  & 11  & 12 \\
6  & $xz$    & $3^{31}$   & N & 3  & 13  & 14  & 15 \\
7  & $xa$    & $3^{0531}$ & P & 4  & 14  & 2  & 15 \\
8  & $xb$    & $0^{31}$   & N & 5  & 15  & 15  &  16\\
9  & $z^2$   & $0^{02}$   & P & 13  & 17  & 17  & 5 \\
10  & $za$   & $0^{420}$  & N & 14  & 17  & 3  & 5 \\
11  & $zb$   & $3^{02}$   & P & 15  & 5  & 5  & 18 \\
12  & $b^2$  & $0^{02}$   & P & 16  & 18  & 18  & 16 \\
13  & $xz^2$ & $1^{13}$   & N & 9  & 19  & 19  &  8\\
14  & $xza$  & $1^{531}$  & N & 10  & 19  & 6  & 8 \\
15  & $xzb$  & $2^{13}$   & N & 11  & 8  & 8  &  20\\
16  & $xb^2$ & $1^{13}$   & N & 12  & 20  & 20  & 12 \\
17  & $z^3$  & $2^{20}$   & N & 19  & 9  & 9  &  11 \\
18  & $zb^2$ & $2^{20}$   & N & 20  & 12  & 12  & 20 \\
19  & $xz^3$ & $3^{31}$   & N & 17  & 13  & 13  & 15 \\
20  & $xzb^2$& $3^{31}$   & N & 18  & 16  & 16  &  18 \\
\end{tabular}
\end{center}
\caption{\label{0123action} Elements, genera, outcomes, and the action of 
generators in the misere quotient semigroup $\mathcal{Q}_{\bf 0.123}$.}
\end{figure}

The outcome class of a given ${\bf 0.123}$ misere position can be determined in two steps.
In the first step, semigroup equivalents for the various single heaps of the position
are looked up using Figure \ref{0123-semi}, whose second row repeats indefinitely.  In the second step,
the semigroup equivalents are multiplied together and the semigroup presentation relations 
in Theorem 1 are used to compute the outcome class.

Here's an example. Suppose we have a sum of six heaps of sizes 
\[1, 3, 4, 8, 9, \mbox{ and } 21.\]  Looking up these heap sizes in 
Figure~\ref{0123-semi}, we find they are equivalent to semigroup elements 
\[x, z, z, a, b, \mbox{ and } b^2,\] respectively.  Multiplying them together, applying commutativity,
and reducing via the semigroup relations, we obtain the element

\[x z^2 a b^3 = (abz)xzb^2 = (b)xzb^2 = (b^3x)z = b^2z = zb^2. \]

Figure \ref{0123action} reveals $zb^2$ to be an N position.

\subsection{\normalsize  Exercise} 
A winning move from the position considered in the previous section 
happens to be to take its entire heap of size 3.  What is the semigroup element and genus of the resulting P position? 
(Answer in this footnote\footnote{[Answer: $b^2$, of genus $0^{02}$.]}).

\section{\normalsize  Knuth-Bendix rewriting and R\'edei's theorem}
The reader might have wondered whether the ``word reduction'' 
algebra of the previous section can always be carried out in general.
Perhaps it's one of those undecidable word problems
that we hear about in semigroup theory?  The answer is no, at least for the particular
example considered---$\mathcal{Q}_{\bf 0.123}$ is a finite semigroup, after all.  
We might as well have written out its entire $20 \times 20$ multiplication table, and used
that to reduce a general word to one of the twenty canonical semigroup elements, 
instead.  

In fact, there's no problem even in the general case, {\em provided} our starting point
is a finitely presented commutative semigroup $S$ 
that we've {\em already proved} to be isomorphic to the misere quotient.  
In this case, there will be no problem with the word reduction problem in $S$ when
we come to apply our semigroup presentation to the selection of best moves in the game.

How does this work in practice?
A {\em canonical presentation} \cite{aks} can be always be computed
from a finitely presented commutative semigroup via the 
{\em Knuth-Bendix rewriting process} \cite{kb}.
In the case of commutative semigroups, the Knuth-Bendix algorithm is always guaranteed to terminate,
and the output of the Knuth-Bendix rewriting process---ie, the canonical presentation---determines an algorithm to solve the word
problem in $S$.

\subsection{\normalsize Confluent rewriting for {\bf 0.123}}

For example, Figure \ref{0123confluent} shows a
{\em confluent rewriting system}
\cite{bn} for ${\bf 0.123}$.  It was computed using
the computer algebra package GAP4 \cite{GAP4}.

\begin{figure}[h]
{\small
\begin{verbatim}
       Rewriting System for Semigroup( [ x, z, a, b, e ] ) with rules 
          [ [ x^2, e ],  [ a^2, e ], , [ z^4, z^2 ], [ a*b, z*b ], 
            [ z^2*b, b ], [ b^3, x*b^2 ], [ z^2*a, z^3 ] ]
\end{verbatim}}
\caption{\label{0123confluent}A confluent rewriting system for ${\mathcal Q}_{\bf 0.123}.$}
\end{figure}

To reduce a general word in the generators $\{x,z,a,b\}$ to an element of 
 ${\mathcal Q}_{\bf 0.123}$, one repeatedly replaces any matching left hand side of a rule
by the corresponding right hand side, until no more such reductions are possible.  
The Knuth-Bendix confluence property guarantees that this process always terminates,
and always terminates
 in the same outcome,
no matter what order the rules are applied.

\subsection{\normalsize R\'edei's theorem}

Another relevant result, due to L. R\'edei, is the following:

{\bf Theorem}
\cite{redei} Every finitely generated commutative semigroup with an identity 
is finitely presentable.

R\'edei's result implies that a {\em partial analysis} of a misere game
(ie, a misere quotient taken only over 
positions with no heap larger than a fixed size $n$) 
will have a finitely presented misere quotient.

\section{\normalsize Pretending}

The reader may have felt another misgiving about our exposition in section \ref{howtolose}.
We started off this paper by highlighting the misere quotient $\mathcal{Q}(\Gamma)$
as the fundamental object of interest in the misere play of a taking and breaking game $\Gamma$.  
But when we actually started to compute outcomes for a concrete, specific position
of {\bf 0.123}, we immediately
bought to bear two additional pieces of information:

\begin{enumerate}
\item Figure \ref{0123-semi} was used to look up a specific semigroup element of
$\mathcal{Q}_{\bf 0.123}$ for each heap in the position; and
\item Figure \ref{0123action}'s partition of $\mathcal{Q}_{\bf 0.123}$ 
elements into P- and N- positions was used to look up the outcome classes of specific
semigroup elements.
\end{enumerate}

Both of these pieces of information are indeed critical to the complete 
analysis of a misere game.
We chose to introduce them informally, first, to simplify our exposition
of the misere quotient machinery.
Now is the time to be more precise.

We call the former information a {\em pretending function} and the latter,
a {\em quotient partition}.  Their definitions 
both involve the misere quotient semigroup $S$.

\subsection{\normalsize Pretending functions \& outcome partitions}
Let $S$ be a semigroup. Let $H$ be the heap alphabet $\{h_1, h_2, \ldots \}$.
A {\em pretending function} is a mapping 
\[\Phi: H \rightarrow S.\]  If $p$ and $r$ are positive integers and $\Phi$ additionally satisfies
\[\Phi(h_k) = \Phi(h_{k+p})\]
for every $k \geq r$, we call $\Phi$ a {\em periodic pretending function of index $r$ period $p$}.

Pretending functions play a role in the analysis of misere games similar
 to the one that {\em nim sequences} do in normal play.

An {\em outcome partition} $S = P \bigcup N$ is 
a partition of a semigroup $S$ into two nonempty
parts, the P positions and the N positions.

\vspace{0.1in}
\begin{figure}[h]
\begin{center}
\begin{tabular}{ll}
Normal play & Misere play \\ \hline
Nim heap equivalent & Quotient semigroup element \\
Nim addition & Quotient semigroup multiplication \\
Periodic nim sequence & Periodic pretending function \\
P position  & P portion of outcome partition \\
\end{tabular}
\end{center}
\caption{Normal vs misere play concepts in impartial games}
\end{figure}

\section{\normalsize A look inside the structure of a wild misere game}
\label{advanced}

In 1940, Rees \cite{r1} proved a fundamental structural result in semigroup theory
that is analogous to the Jordan-H\"older-Schreier Theorem in group theory.  
Rees's theorem asserts that {\em any two relative ideal series of a semigroup $S$ have isomorphic
refinements; in particular, any two composition series of $S$ are isomorphic}\footnote{See \cite{cp}, pg 74.}.  
We will not require (or even state) Rees's results in their full generality, 
but will use ideas associated with Rees's Theorem
to describe the mathematical structure of $\mathcal{Q}_{\bf 0.123}$ in section \ref{0123series}.  Our source 
for the definitions and results cited in this section is \cite{cp}, which we follow closely.

\subsection{\normalsize The Rees congruence.}  Let $I$ be an ideal of a semigroup $S$.  Define a relation $u \ \eta \ v$
to mean that either $u = v$ or else both $u$ and $v$ belong to $I$.  We call $\eta$ the {\em Rees
congruence modulo $I$}.  The equivalence classes of $S \ \mbox{mod } \eta$ are $I$ itself and every one
element set $\{w\}$ with $w$ in $S \setminus I$.  The set $S / I$ should be thought of as the 
result of collapsing $I$ into a single (zero) element, while the elements outside of $I$ retain their identity.

\subsection{\normalsize Principal series \& factors}
\label{principal}
By a {\em principal series} of a semigroup $S$, we mean a chain
\begin{equation}
\label{0123principal}
S = S_1 \supset S_2 \supset \cdots \supset S_m \supset S_{m+1} = \emptyset 
\end{equation} 
of ideals $S_i$ $(i = 1, \ldots, m)$ of $S$, beginning with $S$ and ending with the empty set, and such
that there is no ideal of $S$ strictly between $S_i$ and $S_{i+1}$ $(i=1,\ldots,m)$.  

By the {\em factors} of a
principal series $S$, we mean the Rees factor semigroups $S_i/S_{i+1}$.

\subsection{\normalsize Principal series for $\mathcal{Q}_{\bf 0.123}$}
\label{0123series}

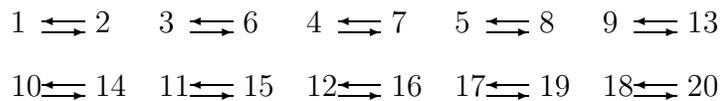
\begin{figure}[h]
\begin{center}
\begin{picture}(45,40)(0,10)
\setlength{\unitlength}{.8\unitlength}

\put(-135,10){10}
\put(-120,10){\vector(1,0){20}}
\put(-100,15){\vector(-1,0){20}}
\put(-95,10){14}

\put(-65,10){11}
\put(-50,10){\vector(1,0){20}}
\put(-30,15){\vector(-1,0){20}}
\put(-25,10){15}

\put(5,10){12}
\put(20,10){\vector(1,0){20}}
\put(40,15){\vector(-1,0){20}}
\put(45,10){16}

\put(75,10){17}
\put(90,10){\vector(1,0){20}}
\put(110,15){\vector(-1,0){20}}
\put(115,10){19}

\put(145,10){18}
\put(160,10){\vector(1,0){20}}
\put(180,15){\vector(-1,0){20}}
\put(185,10){20}

\put(-135,40){1}
\put(-120,40){\vector(1,0){20}}
\put(-100,45){\vector(-1,0){20}}
\put(-95,40){2}

\put(-65,40){3}
\put(-50,40){\vector(1,0){20}}
\put(-30,45){\vector(-1,0){20}}
\put(-25,40){6}

\put(5,40){4}
\put(20,40){\vector(1,0){20}}
\put(40,45){\vector(-1,0){20}}
\put(45,40){7}

\put(75,40){5}
\put(90,40){\vector(1,0){20}}
\put(110,45){\vector(-1,0){20}}
\put(115,40){8}

\put(145,40){9}
\put(160,40){\vector(1,0){20}}
\put(180,45){\vector(-1,0){20}}
\put(185,40){13}

\end{picture}
\end{center}
\caption{\label{0123x} The action of $x$ on $\mathcal{Q}_{\bf 0.123}$.}
\end{figure}

The final four columns of Figure \ref{0123action} show the action of each of the four generators $\{x,z,a,b\}$ on
$\mathcal{Q}_{\bf 0.123}$.  We've written the same information graphically in Figures
\ref{0123x}, \ref{0123z}, \ref{0123a}, and \ref{0123b}.   The elements $x$ and $a$ are simple involutions,
but $z$ and $b$ are not.


\begin{figure}[h]
\begin{center}
\begin{picture}(45,60)(0,10)
\setlength{\unitlength}{.8\unitlength}

\put(-75,70){1}
\put(-60,72){\vector(1,0){20}}

\put(-35,70){3}
\put(-20,72){\vector(1,0){20}}

\put(5,70){9}
\put(20,70){\vector(1,0){20}}
\put(40,75){\vector(-1,0){20}}
\put(45,70){17}

\put(85,70){10}
\put(80,72){\vector(-1,0){20}}

\put(125,70){4}
\put(120,72){\vector(-1,0){20}}

\put(-75,40){2}
\put(-60,42){\vector(1,0){20}}

\put(-35,40){6}
\put(-20,42){\vector(1,0){20}}

\put(5,40){13}
\put(20,40){\vector(1,0){20}}
\put(40,45){\vector(-1,0){20}}
\put(45,40){19}

\put(85,40){14}
\put(80,42){\vector(-1,0){20}}

\put(125,40){7}
\put(120,42){\vector(-1,0){20}}

\put(-115,10){5}
\put(-100,10){\vector(1,0){20}}
\put(-80,15){\vector(-1,0){20}}
\put(-75,10){11}

\put(-45,10){8}
\put(-30,10){\vector(1,0){20}}
\put(-10,15){\vector(-1,0){20}}
\put(-5,10){15}

\put(25,10){16}
\put(40,10){\vector(1,0){20}}
\put(60,15){\vector(-1,0){20}}
\put(65,10){20}

\put(95,10){12}
\put(110,10){\vector(1,0){20}}
\put(130,15){\vector(-1,0){20}}
\put(135,10){18}

\end{picture}
\end{center}
\caption{\label{0123z} The action of $z$ on $\mathcal{Q}_{\bf 0.123}$.}
\end{figure}


\begin{figure}[h]
\begin{center}
\begin{picture}(45,40)(0,10)
\setlength{\unitlength}{.8\unitlength}

\put(-135,10){5}
\put(-120,10){\vector(1,0){20}}
\put(-100,15){\vector(-1,0){20}}
\put(-95,10){11}

\put(-65,10){6}
\put(-50,10){\vector(1,0){20}}
\put(-30,15){\vector(-1,0){20}}
\put(-25,10){14}

\put(5,10){9}
\put(20,10){\vector(1,0){20}}
\put(40,15){\vector(-1,0){20}}
\put(45,10){17}

\put(75,10){13}
\put(90,10){\vector(1,0){20}}
\put(110,15){\vector(-1,0){20}}
\put(115,10){19}

\put(145,10){16}
\put(160,10){\vector(1,0){20}}
\put(180,15){\vector(-1,0){20}}
\put(185,10){20}

\put(-135,40){1}
\put(-120,40){\vector(1,0){20}}
\put(-100,45){\vector(-1,0){20}}
\put(-95,40){4}

\put(-65,40){2}
\put(-50,40){\vector(1,0){20}}
\put(-30,45){\vector(-1,0){20}}
\put(-25,40){7}

\put(5,40){3}
\put(20,40){\vector(1,0){20}}
\put(40,45){\vector(-1,0){20}}
\put(45,40){10}

\put(75,40){8}
\put(90,40){\vector(1,0){20}}
\put(110,45){\vector(-1,0){20}}
\put(115,40){15}

\put(145,40){12}
\put(160,40){\vector(1,0){20}}
\put(180,45){\vector(-1,0){20}}
\put(185,40){18}

\end{picture}
\end{center}
\caption{\label{0123a} The action of $a$ on $\mathcal{Q}_{\bf 0.123}$.}
\end{figure}


\begin{figure}[h]
\begin{center}
\begin{picture}(45,150)(0,-70)
\setlength{\unitlength}{.8\unitlength}


\put(-37,74){10}
\put(-32,72){\vector(0,-1){20}}

\put(82,74){13}
\put(88,72){\vector(0,-1){20}}

\put(-75,40){1}
\put(-60,42){\vector(1,0){20}}

\put(-35,40){5}
\put(-20,42){\vector(1,0){20}}

\put(5,40){12}
\put(20,40){\vector(1,0){20}}
\put(40,45){\vector(-1,0){20}}
\put(45,40){16}

\put(85,40){8}
\put(80,42){\vector(-1,0){20}}

\put(125,40){2}
\put(120,42){\vector(-1,0){20}}

\put(-35,6){9}
\put(-32,17){\vector(0,1){20}}

\put(82,6){14}
\put(87,17){\vector(0,1){20}}

\put(-37,-20){17}
\put(-32,-22){\vector(0,-1){20}}

\put(88,-20){7}
\put(91,-22){\vector(0,-1){20}}

\put(-75,-60){3}
\put(-60,-58){\vector(1,0){20}}

\put(-35,-60){11}
\put(-20,-58){\vector(1,0){20}}

\put(5,-60){18}
\put(20,-60){\vector(1,0){20}}
\put(40,-55){\vector(-1,0){20}}
\put(45,-60){20}

\put(85,-60){15}
\put(80,-58){\vector(-1,0){20}}

\put(125,-60){6}
\put(120,-58){\vector(-1,0){20}}

\put(-35,-96){4}
\put(-32,-86){\vector(0,1){20}}

\put(86,-96){19}
\put(91,-86){\vector(0,1){20}}
\end{picture}
\end{center}
\caption{\label{0123b} The action of $b$ on $\mathcal{Q}_{\bf 0.123}$.}
\end{figure}


Such pictures aid in computing the transformations
corresponding to the other sixteen semigroup elements of $\mathcal{Q}_{\bf 0.123}$, 
finding the principal ideal series of the semigroup, and working out the associated Rees factor semigroups.

It's not too hard to show (eg, \cite{lalle}, Proposition 1.3, pg 21) that a semigroup $S$
of transformations of a finite set $X$ (in this example, $X=\{1,2,\cdots,20\}$) has a unique minimal ideal $J$, and that $J$ necessarily
consists of the elements $u \in S$ with the smallest cardinality image set (or {\em rank})  $|uS|$.
For $S=\mathcal{Q}_{\bf 0.123}$, these smallest-rank elements are
$J= \{12, 16, 18, 20\}$. They each have rank four.
The minimal (or {\em kernel}) ideal is therefore
\[J= S_5 = \{b^2,xb^2,zb^2,xzb^2\},\]
which is part of the entire principal series 
\[\mathcal{Q}_{\bf 0.123} = S = S_1 \supset S_2 \supset S_3 \supset S_4 \supset S_5 \supset S_{6} = \emptyset,\]
where
\begin{eqnarray*}
S_1 & = & S_2 \cup \{e, x, a, xa\}. \\
S_2 & = & S_3 \cup \{z, xz, za, xza \} \\
S_3 & = & S_4 \cup \{z^2,xz^2,z^3,xz^3\} \\
S_4 & = & S_5 \cup \{b,xb,zb,xzb\} \\
S_5 & = & \{b^2,xb^2,zb^2,xzb^2\} \\
S_6 & = & \emptyset. \\
\end{eqnarray*} 

Given such a principal series, we'll also define $D_n = S_n \setminus S_{n+1}$.  The
sets $D_n$ partition $S$, and together form the congruence classes of
the {\em mutual divisibility congruence} on $S$ (see \S \ref{divisibility}).

The Rees factor semigroups $S_1/S_2$, $S_3/S_4$, and $S_5/S_6$ are each isomorphic to 
the {\em Klein four-group} ${Z}_2 \times {Z}_2$ {\em with adjoined zero}
\[K_4 \cup \{0\} = \left({Z}_2 \times {Z}_2 \right) \cup \{0\}.\] 
For example, Figure \ref{zmult} contains the multiplication table of $S_3/S_4$.

\begin{figure}[h]
\begin{center}
\begin{tabular}{c|ccccc}
     & 0 & $z^2$ & $xz^2$ & $z^3$ & $xz^3$ \\ \hline
0    & 0 &  0  &  0   &  0  &  0   \\  
$z^2$  & 0 &  $z^2$ &  $xz^2$ & $z^3$  & $xz^3$   \\  
$xz^2$ & 0 & $xz^2$  &  $z^2$   &  $xz^3$  &  $z^3$   \\  
$z^3$  & 0 & $z^3$  & $xz^3$   &  $z^2$  &  $xz^2$   \\  
$xz^3$ & 0 &  $xz^3$  &  $z^3$   &  $xz^2$  &  $z^2$   \\
\end{tabular}
\end{center}
\caption{\label{zmult} Multiplication in $S_3/S_4$.}
\end{figure}

The group $K_4$ also happens to be the {\em normal play} quotient of ${\bf 0.123}$.  One 
often finds such an isomorphism between the normal play quotient of a game and one of its misere game Rees factors (after
deleting the adjoined zero)---in fact, illustrating this phenomenon 
was our motivation for calculating the Rees factors in the first place!

The factors $S_2/S_3$ and $S_4/S_5$ are {\em null semigroups}---all products equal zero.
See Figure \ref{zamult}.

\begin{figure}[h]
\begin{center}
\begin{tabular}{c|ccccc}
     & 0 & $b$ & $xb$ & $zb$ & $xzb$ \\ \hline
0    & 0 &  0  &  0   &  0  &  0   \\  
$b$  & 0 &  0  &  0   &  0  &  0   \\  
$xb$ & 0 &  0  &  0   &  0  &  0   \\  
$zb$  & 0 &  0  &  0   &  0  &  0   \\  
$xzb$ & 0 &  0  &  0   &  0  &  0   \\  
\end{tabular}
\end{center}
\caption{\label{zamult} Multiplication in $S_4/S_5$.}
\end{figure}

\section{\normalsize Idempotents \& tame islands}

Using some more elementary definitions and semigroup results,  
it's possible to 
shed more light on the structure of tame and wild positions in a misere quotient such as 
${\mathcal Q}_{\bf 0.123}$.

\subsection{\normalsize The natural partial ordering of idempotents}
A reflexive, antisymmetric, and transitive 
relation $\leq$ on a set $X$ is called a {\em partial ordering}. 
An element $f$ of a semigroup $S$ is an {\em idempotent} if $f^2 = f$.

Let $E$ be the set of all idempotents of a semigroup $S$.  Define $g \leq f$ (for $g,f \in E$) to mean
$gf = fg = g$.  Then $\leq$ is a partial ordering of $E$ which we call the {\em natural partial ordering
of idempotents of $S$}.  (See \cite{cp}, its \S 1.8, for a proof that $\leq$ is a partial ordering of $E$).

\vspace{0.1in}
Let's apply these definitions to $\mathcal{Q}_{\bf 0.123}$.  After computing the square of each its twenty elements
(Figure \ref{0123action}) to see if it is an idempotent,
one finds that $E = \{e, z^2, b^2\}$.  It's easy to see that 
\[ z^2 \leq e \mbox{, since } z^2 e = z^2,\]
and
\[ b^2 \leq e \mbox{, since } b^2 e = b^2.\]
How do $z^2$ and $b^2$ compare?  Start with $z^2b^2$, insert
an $a^2 = e$ factor, and apply the semigroup relation $zab=b$ (\S \ref{0123}) twice:

\begin{eqnarray*}
z^2b^2 & = & z^2(a^2)b^2 \\
       & = & (zab)(zab) \\
       & = & b^2,
\end{eqnarray*}
ie,
\[b^2 \leq z^2.\]
The natural 
partial ordering
of idempotents of $\mathcal{Q}_{\bf 0.123}$ is therefore the linear ordering  
\[b^2 \leq z^2 \leq e.\]

\subsection{\normalsize Divisibility \& idempotents: the tame islands.}

\label{divisibility}
If $u$  and $v$ are elements of a commutative semigroup $S$ we say that $u$ {\em divides} $v$, and write $u | v$, 
if there exists a $w$ such that $uw=v$.  Define a relation $u \ \tau \ v$ to mean that both $u$ divides $v$ 
and $v$ divides $u$.  We'll call $\tau$ the {\em mutual divisibility} relation.

\vspace{0.1in}

Here's a useful result (\cite{cp}, pg 22):

{\bf Theorem}
A semigroup $S$ contains a subgroup if and only if it contains an idempotent.

\vspace{0.1in} 

Here's another useful result \cite{aks}:

{\bf Theorem}
In a commutative semigroup $S$, the 
mutual divisibility relation is a congruence.
The congruence class $f \tau$ containing an idempotent $f$
is precisely the maximal subgroup of $S$ for which that element is the 
identity.

Applying the latter theorem to the principal series of $\mathcal{Q}_{\bf 0.123}$ 
(equation \ref{0123principal}),
we find that the sets 
\begin{eqnarray*}
D_1 & = & \{e, x, a, xa\}. \\
D_3 & = & \{z^2,xz^2,z^3,xz^3\} \\
D_5 & = & \{b^2,xb^2,zb^2,xzb^2\} \\
\end{eqnarray*} 
are three disjoint {\em subgroups} of $\mathcal{Q}_{\bf 0.123}$.  Their respective
identities 
\[e, z^2, b^2\]  
are the three idempotents of $\mathcal{Q}_{\bf 0.123}$.
The semigroup
multiplication for elements within each of the subgroups $D_1$, $D_3$, and $D_5$
follows that of 
${Z}_2 \times {Z}_2$,
which in turn is just the same as that of the misere Nim positions 
\[\{\ast 2+ \ast 2, \ast 1, \ast 2, \ast 3 \}\]
with respective genera
\[\{0^{02}, 1^{13}, 2^{20}, 3^{30} \}.\]
We call $D_1$, $D_3$ and $D_5$ the {\em tame islands} of $\mathcal{Q}_{\bf 0.123}$.

\subsection{\normalsize Decomposition of a commutative semigroup}

The type of commutative semigroup decomposition carried out in the previous
section was first completely described in 1954 
by Tamura and Kimura \cite{tk}.  To describe their
result, we need one more definition:

A semigroup $S$ is {\em archimedean} if for any two elements of $S$, each divides 
some power of the other.

Here is the result of Tamura and Kimura (see also \cite{cp}, its \S 4.3, pg 135):

{\bf Theorem}
 Any commutative semigroup $S$ is uniquely expressible as a semilattice $Y$ of
 archimedean semigroups $S_\alpha$.  The semigroup $S$ can be embedded in a 
 semigroup $T$ which is a union of groups if and only if $S$ is separative,
 and this is so if and only if each $S_\alpha$ is cancellative.  The semigroup
 $T$ can be taken to be the union of the same semilattice $Y$ of groups $G_\alpha$,
 where $G_\alpha$ is the quotient group of $S_\alpha$, for each $\alpha$.

We needn't concern ourselves too much with the 
theorem of Tamura and Kimura, except to the extent that it 
informs us of the most general structure conceivable for a misere
game's quotient.   Our interests lie in computing 
quotient semigroup presentations for particular
games, casting their
solutions in the form of periodic pretending functions. 
In the general case, this may be difficult---just as it is in many normal play games.  
Nevertheless, many weapons at hand in normal play have their natural analogues in misere play.
The remainder of this paper is therefore devoted to computational
aspects of quotient semigroup construction.  

\section{\normalsize Proving quotients correct}
\label{proving}
In 1955, Guy and Smith
proved a result about normal play {\em octal games} (\cite{ww1}, Chapter 4)
that they stated as follows:

{\bf Theorem}  [\cite{gs}, pg 516]:
If a game $\Gamma$ is defined by a finite octal, having $P$
places after the point, and if we can empirically find positive integers
$p$ and $r_0$ such that the equation
\[G(r+p) = G(r)\] is true for all $r$ in the range $r_0 \leq r < 2r_0+p+P$,
then it is true for all $r \geq r_0$, so that $G$-function has ultimate
period $p$.

Our goal in this section is prove an analogue of the Guy \& Smith result for misere play.
Roughly speaking, this involves
replacing the normal-play single-heap nim 
equivalence function $G()$ by an appropriate periodic pretending function $\Phi()$,
and replacing the notion of normal play P-positions corresponding to positions of nim equivalent $\ast 0$ with
an appropriate outcome partition of ${\mathcal Q}(\Gamma)$.

\subsection{\normalsize Correctness to heap size $n$}

Fix a positive integer $n$ and heap alphabet $H_n = \{h_1, h_2, \ldots, h_n\}$.  
We first consider the problem of verifying
that an asserted {\em finite} and {\em explicit} 
misere quotient semigroup ${\mathcal Q}(\Gamma)$, pretending function $\Phi$, and
outcome partition are {\em correct to heap size $n$} in the sense that they together
correctly describe all P- and N-positions
of the misere play of a game $\Gamma$, {\em provided no heap is larger than size $n$}.  

Here's the necessary machinery.

\subsubsection{\normalsize Move pairs}
Suppose 
\begin{equation}
\label{genmove}
 h_f \rightarrow t 
\end{equation}
is a concrete move of $\Gamma$ that involves replacing the single heap of
size $f \leq n$ with various other smaller heaps represented by the element $t \in {\mathcal F}_{H_n}$, according to the rules.  The given 
pretending
function $\Phi$ determines a pair of corresponding semigroup elements
$(s_f, s)$ in the obvious way: on the left-hand side of (\ref{genmove}), 
$s_f$ is just the pretending function image of 
the heap $h_f$; on the right-hand side, we multiply the images of the various heaps $h_i$ occuring in $t$ together to obtain $s$.
Each $(s_f, s) \in Q \times Q$ that can be formed in this way
will be called a {\em move pair to heap size $n$}.  

The set of $M_n$ of all move pairs to heap size $n$ can be thought of as a relation on ${\mathcal Q}$, but 
we're {\em not} expecting it to be a reflexive or symmetric relation, since
each pair is derived from a move in the game, and these each have a ``direction.''

{\bf Lemma}
The set $M_n$ of all move pairs to heap size $n$ is finite.

{\it Proof} \/
Since there are only finitely many moves to heap size $n$, the result is immediate.

\subsubsection{\normalsize Move pair translates}
What want to do with $M_n$ is  extend it to all pairs of the form
\[ (u \cdot s_f, u \cdot s), \]
where $u$ is an arbitrary element of $\mathcal{Q}$. 
This is the set $T_n$
of all {\em move pair translates to heap size $n$}
\[ T_n = \bigcup_{\stackrel{{\textstyle u \in \mathcal{Q}}}
                         {{\textstyle (s_f,s) \in M_n}}} 
                         (u \cdot s_f, u \cdot s), \]
where $f$ is allowed to range freely $1 \leq f \leq n$.
We sometimes call the quotient semigroup element $u$ the {\em basis} of a translate.  

{\bf Lemma}
The set $T_n$ of all move translates is a finite relation on ${\mathcal Q}$.

{\it Proof} \/
Since there are only finitely many elements in ${\mathcal Q}$ by assumption, the result is immediate.

\subsubsection{\normalsize Verification algorithms}
\label{verification}
Now let's consider what we need to do in order to prove that an analysis is
correct to heap size $n$.  An induction argument needs two subpieces to
be successful:

\begin{enumerate}
\item Show that there is no move from a concrete position asserted to be P to
another P position.
\item Show that every non-endgame position asserted to be N has some move to
a P position.
\end{enumerate}

We can dispense with the first case
by computing all the move pair translates to heap size $n$, and seeing
whether there is any translate of the form 

\[(\mbox{P position}, \mbox{P position}).\] 

Confirming that there is no such translate already completes the first
half the inductive argument.  An example computation of this type is carried out
in \S \ref{PtoP}, below.

For the second half, we want to make sure that we can always find
a move from every non-endgame position asserted to be an $N$-position to some $P$-position.
This is more complicated.  In the algorithm to be described, each
$N$-position type $\omega$ in the given outcome partition of 
${\mathcal Q}$ is considered separately.  Roughly speaking, the desired algorithm works by
examining  of the ``fine structure'' of move pair translates of the form $(\omega,P)$,
restricted to each subset $U$ of $H_n$.   

In order to be precise, we need some more definitions
and notation.

For a given nonempty set $U$ with $U \subseteq H_n$, we define

\[B(U) = \{p \in {\mathcal F}_{U} \ | \ \forall \ h_i \in U, \mbox{ the heap } h_i \mbox{ occurs at least once in } p\}.\]

There are $2^n-1$ such sets $B(U)$, and together they form a partition of the infinite set ${\mathcal F}_{H_n}$ into
a finite number of parts.  We also define a symbol for the product of the elements of $U$:
\[P(U) = \prod_{h_i \in U} h_i.\]

If $h_i \in U$, we define an operator $\partial / \partial h_i$ acting on positions $p \in B(U)$ by
\[\frac{\partial}{\partial h_i} p = \hat{p},\]
where $\hat{p} = p / h_i$, ie, $\hat{p}$ is the position obtained by deleting one heap of size $i$ from $p$.
It follows that for $p \in B(U)$ and each $h_i \in U$,
\[p = h_i \frac{\partial}{\partial h_i} p.\]

A general element $p \in B(U)$ can always be written (up to the ordering of factors) in the form
\begin{equation}
\label{rprod}
p = r \prod_{h_i \in U} h_i = r \cdot P(U)
\end{equation}
for a unique value $r \in {\mathcal F}_{H_n}$.  (It may be that $r = \Lambda = 1$ is the empty position, i.e., the endgame).
In fact, if $U$ is the set of heaps of sizes $\{i_1,\cdots, i_k\}$, then
\[ r = \frac{\partial}{\partial h_{i_1}}\frac{\partial}{\partial h_{i_2}}\cdots\frac{\partial}{\partial h_{i_k}} \ p.\]

A general legal move from a position $p \in B(U)$ is of the form
\begin{equation}
\label{genmoveform}
h_i \rightarrow t,
\end{equation}
where $h_i \in U$ and the allowable values $t \in {\mathcal F}_{H_n}$   
are determined by the rules $\Gamma$.  Define a corresponding 
subset $M_n(U) \subseteq M_n$ of move pairs of $\Gamma$
\begin{equation}
\label{genmovepair}
M_n(U) = \{\ (\Phi(h_i),\Phi(t)) \ | \ h_i \in U, \mbox{ and } h_i \rightarrow t \mbox{ is a legal move in the play of } \Gamma\}.
\end{equation}

Let ${\mathcal S}(U)$ be the subsemigroup of ${\mathcal Q}$ generated by the identity $e$ of $\mathcal Q$ and all the
$\Phi(h_i)$'s for $h_i \in U$.
We're going to be interested in a set of move pair translates $T_n(U) \subseteq T_n$ obtained from the set described
in equation (\ref{genmovepair})

\[T_n(U) = \{\ (s \cdot \Phi(h_i),s \cdot \Phi(t)) \ | \ (\Phi(h_i),\Phi(t)) \in M_n(U), \mbox{ and } s \in {\mathcal S}(U)\}.\]

Finally, abusing our notation slightly for notational convenience, for a nonempty $U \subseteq H_n$, we define
\[\Phi(U) = \prod_{h_i \in U} \Phi(h_i) = \Phi(P(U)), \] and define $\Phi(U)=1$ if $U$ is empty.

We're ready to prove the following theorem.

{\bf Theorem 2}.  Suppose $\omega \in {\mathcal Q}$ is asserted to be an N-position type.  The following are equivalent.
\begin{enumerate}
\item Every non-endgame position $p \in {\mathcal F}_{H_n}$ asserted to be of  type $\omega$ has a move to a P-position.
\item For every choice of a nonempty subset $U \subseteq H_n$ and quotient semigroup element $s \in {\mathcal S(U)}$ such that
the equation
\begin{equation}
\label{divispossibility}
\omega = s \cdot \Phi(U)
\end{equation} is satisfied, there is some heap $h_i \in U$ and legal move of $\Gamma$ 
\[h_i \rightarrow t\]
such that
\begin{equation}
\label{pwinning}
(s \cdot \Phi(U)\ ,\ s \cdot \Phi(t) \cdot \Phi(\frac{\partial}{\partial h_i} \ P(U))
\end{equation}
is a move pair translate of the form
\[(\omega, P\mbox{-position}).\]
\end{enumerate}

{\it Proof.}\/
$(1) \Rightarrow (2)$.  
Let $U = \{h_{i_1},\cdots,h_{i_k}\}$ and $s \in {\mathcal S}(U)$ 
be an arbitrary
solution of equation (\ref {divispossibility}) for some $k \geq 1$.  By definition
we can write
\[s = \prod_{1 \leq j \leq k} \Phi(h_{i_j})^{\nu_j}\]
for some appropriate powers $\nu_j \geq 0$.  Forming a position $p \in {\mathcal F}_{H_n}$ as the product
\begin{equation}
\label{explicitpos}
p = \left( \prod_{1 \leq j \leq k} h_{i_j}^{\nu_j} \right) \left( \prod_{1 \leq j \leq k} h_{i_j} \right),
\end{equation}
we have $\Phi(p) = s \cdot \Phi(U) = \omega$.  If $p$ is not the endgame, by assumption there is a winning move 
\begin{equation}
\label{wm}
h_{i_j} \rightarrow t,
\end{equation}
and we observe that this move can be made in the
 second half of the product on the right
hand side of (\ref{explicitpos}).  The resulting
$P$-position $q$ is of the form  
\begin{equation}
\label{formofp}
q = \left( \prod_{1 \leq j \leq k} h_{i_j}^{\nu_j} \right) \left( \frac{\partial}{\partial h_{i_j}} \prod_{1 \leq j \leq k} h_{i_j} \right) \ t. 
\end{equation}

In particular, $(\Phi(h_{i_j}),\Phi(t))$ is a move pair to heap size $n$,
and translating that move pair by the basis
\[b = s \cdot \Phi(\frac{\partial}{\partial h_i} \ P(U))\]
yields a move pair translate
\[(s \cdot \Phi(U)\ ,\ s \cdot \Phi(t) \cdot \Phi(\frac{\partial}{\partial h_i} \ P(U))\]
of the form
\[(\omega, P\mbox{-position}),\]
as desired.

$(2) \Rightarrow (1)$.   The converse is similar.  Suppose $p \in {\mathcal F_{H_n}}$ is a non-endgame position of type $\omega$.
Let \[U = \{ h_i \ | \ h_i \mbox{ occurs at least once as a heap in } p \}.\]  Then $p$ can be written
in the form of equation (\ref{rprod}) for a unique $r \in {\mathcal F_{H_n}}$, and applying the pretending
function $\Phi$ to both sides of equation (\ref{rprod}), we obtain equation (\ref{divispossibility}), where $s = \Phi(r)$.
The hypothesis applies, so there is a legal winning move to a P position
\[h_i \rightarrow t\]
such that applying that move to $p$ leads to a position $q$ of the form in equation (\ref{formofp}) with corresponding
move translate of the form in equation (\ref{pwinning}).  So $p$ has a move to a $P$-position, as desired.  \framebox[2mm]

Relying upon Theorem 2, we can now describe a brute-force algorithm for verifying that every non-endgame position asserted to be
an $N$ position has some move to a $P$ position, to heap size $n$.  Simply put, each $N$ position type $\omega$ is considered separately, and
must pass the test given in the second half of the statement of Theorem 2.  Since (i) the number $2^n$ of
subsets $U$ to be considered is finite; and (ii) there's
only a finite number of moves to heap size $n$; and finally, (iii) there are only a finite number of choices for values $s$,
since they're taken from various subsemigroups of the 
proposed quotient semigroup ${\mathcal Q}$,
which required to be finite by assumption,
we can exhaust all possible ways of forming equation (\ref{divispossibility}).  If the desired translates
exist for every value of $\omega$, $U$, $s$, and move $h_i \rightarrow t$ meeting the conditions stated in Theorem 2,
the verification succeeds; otherwise, it fails.

An illustration computation of the type described in Theorem 2 is carried out in section \S \ref{NtoP}, below.

\subsubsection{Complexity}

The algorithm given at the beginning of section \ref{verification} 
for verifying that no $P \rightarrow P$ moves exist to heap size $n$ is a polynomial-time
computation in its natural parameters $n$ and $|{\mathcal Q}|$.  By contrast, the complexity of the 
verification algorithm we've given for $N \rightarrow P$ moves 
is {\em exponential} in $n$, since it involves considering all subsets $U \subseteq H_n$.  In practice, the latter algorithm can be sped up dramatically by exploiting 
available information about the values present in the proposed pretending function $\Phi$ and the structure of 
the quotient semigroup ${\mathcal Q}$.  For example,
the $N \rightarrow P$ verification algorithm to heap size $n$ as presented in this paper ignores all information that may be already
available about the correctness of the proposed quotient to heap size $n-1$.   
A more complete discussion of such improvements, which involve a combination of semigroup theory,
data-structure organization, and heuristics, would take us too far afield from the goals of this paper.

Closely related to the problem of $N \rightarrow P$ move verification is 
the important problem of {\em misere quotient semigroup construction}.    This is another
rich subject with many points of contact with semigroup structure theory and algorithm design, but 
it is not treated in this paper.

\subsection{\normalsize Correctness of the {\bf 0.123} quotient to heap size 12}

Let's illustrate how the two verification computations work in the specific case of ${\bf 0.123}$, to
heap size 12.  It less than one second to run both verification procedures
in GAP \cite{GAP4}, a computer algebra package.  

\subsubsection{\normalsize No $P \rightarrow P$ moves exist to heap size 12}

We include in Figure \ref{PP123} a portion of an automatically generated proof that there's
``no move from a position asserted to P to a P position'' to heap size 12 in ${\bf 0.123}$.
\label{PtoP}

Recall that the positions asserted to be P are the semigroup elements
\[P = \{x, xa, z^2, zb, b^2\}.\]

If there were a move from a position asserted to P to a P position,
there would be a translate $(u \cdot s_f, u \cdot s)$ of the form 
$(\mbox{P-position},\mbox{P-position})$.

Figure 12 lists all moves to heap size twelve, and the associated translates in the particular subcase where the
basis of the translate sets is the semigroup element $x$.   Since there is no translate of the form $(P,P)$,
we are done.  The computation for the other 19 basis elements in the semigroup is similar.

\begin{figure}[h]
\begin{center}
\begin{verbatim}
            Move     Move Pair Translated by x       Outcomes
            1->0             (e,x)                    (N,P)
            3->0             (x*z,x)                  (N,P)
            3->1             (x*z,e)                  (N,N)
            4->1             (x*z,e)                  (N,N)
            4->2             (x*z,x)                  (N,P)
            5->2             (e,x)                    (N,P)
            5->3             (e,x*z)                  (N,N)
            6->3             (x*b^2,x*z)              (N,N)
            6->4             (x*b^2,x*z)              (N,N)
            7->4             (x,x*z)                  (P,N)
            7->5             (x,e)                    (P,N)
            8->5             (x*a,e)                  (P,N)
            8->6             (x*a,x*b^2)              (P,N)
            9->6             (x*b,x*b^2)              (N,N)
            9->7             (x*b,x)                  (N,P)
            10->7            (e,x)                    (N,P)
            10->8            (e,x*a)                  (N,P)
            11->8            (x*b^2,x*a)              (N,P)
            11->9            (x*b^2,x*b)              (N,N)
            12->9            (x,x*b)                  (P,N)
            12->10           (x,e)                    (P,N)
\end{verbatim}
\end{center}
\caption{\label{PP123} Verification that there is no $P \rightarrow P$ translate to heap size 12 in
{\bf 0.123} in the particular case that the basis is equal to $x$.  }
\end{figure}

\subsubsection{\normalsize $N \rightarrow P$ moves must exist, to heap size 12}
\label{NtoP}
In this section, we illustrate how to verify that a position asserted to be $N$ must have
a move to a $P$ position, again using ${\bf 0.123}$ as our example.  Because the algorithm we've
presented in \S\ref{verification} involves considering each of the $2^{12} -1 = 4095$ nonempty subsets $U$
of $H_{12}$ in turn, we must satisfy ourselves by showing only a small portion of the computation.

Take, for example, the semigroup element $\omega = xb$, 
which is asserted be an $N$ position in ${\bf 0.123}$.  Theorem 2 informs us
that we're going to be interested in move translates of the form

\[(xb,P\mbox{-position}).\]

All such move translates are shown in the first column of Figure \ref{NtoPfig}.
Each is listed together with the basis element $u$, concrete move, and move pair
that generates it.
\begin{figure}[h]
\begin{center}
\begin{tabular}{cccc}
 Move         &               &                    &           \\         
Translate & Basis $u$    &     Move           & Move pair \\              
$(xb,x)$  & $x$           & $9 \rightarrow 7$  & $(b,e)$  \\
$(xb,zb)$ & $b$           & $5 \rightarrow 3$  & $(x,z)$ \\
$(xb,zb)$ & $b$           & $10 \rightarrow 8$ & $(x,a)$  \\
$(xb,zb)$ & $xzb$         & $3 \rightarrow 1$  & $(z,x)$  \\
$(xb,zb)$ & $xzb$         & $8 \rightarrow 5$  & $(a,x)$  \\
$(xb,zb)$ & $xzb$         & $4 \rightarrow 1$  & $(z,x)$  
\end{tabular}
\caption{\label{NtoPfig} The translates of the form $(xb,\mbox{P-position})$ 
for ${\bf 0.123}$ with no heap larger than size 12.}
\end{center}
\end{figure} 

To illustrate the main features of the verification computation for positions of type $\omega = xb$,
we'll show operation of the algorithm on the fifteen nonempty subsets $U$ of $\{h_4,h_8,h_{9},h_{10}\}$.

\begin{figure}[h]
\begin{center}
\begin{tabular}{cc}
$U$ & ${\mathcal S}(U)$ \\ \hline
$\{h_4\}$ & $\{1,z,z^2,z^3\}$ \\
$\{h_8\}$ & $\{1,a\}$ \\
$\{h_9\}$ & $\{1,b,b^2,xb^2\}$ \\
$\{h_{10}\}$ & $\{1,x\}$ \\
$\{h_4,h_8\}$ & $\{1,z,a,z^2,az,z^3\}$ \\
$\{h_4,h_9\}$ & $\{1,z, b, z^2, bz, b^2, z^3, b^2z, xb^2, b^2xz \}$ \\
$\{h_4,h_{10}\}$ & $\{1, z, x, z^2, xz, z^3, xz^2, xz^3 \}$ \\
$\{h_8,h_9\}$ & $\{1,a, b, bz, b^2, b^2z, b^2x, b^2xz\}$ \\
$\{h_8,h_{10}\}$ & $\{1,x,a,ax\}$ \\
$\{h_9,h_{10}\}$ & $\{\framebox{1}, b, x, b^2, bx,  b^2x\}$ \\
$\{h_4,h_8,h_9\}$ & $\{1, z, a, b, z^2, az, bz, b^2, z^3, b^2z, b^2x, b^2xz\}$ \\
$\{h_4,h_8,h_{10}\}$ & $\{ 1, z, a, x, z^2, az, xz, ax, z^3, xz^2, axz, xz^3 \}$ \\
$\{h_4,h_9,h_{10}\}$ & $\{1,\framebox{$z$}, b, x, z^2, bz, xz, b^2, bx, z^3, xz^2, b^2z, bxz, b^2x, xz^3, b^2xz\}$ \\
$\{h_8,h_9,h_{10}\}$ & $\{1, b, \framebox{$a$}, x, b^2, bz, bx, ax, b^2x, b^2z, bxz, b^2xz\}$ \\
$\{h_4,h_8,h_9,h_{10}\}$ & $\{\framebox{$1$}, a, b, x, z, bz, ax,\framebox{$az$}, b^2, bx, xz, \framebox{$z^2$},$ \\
                         & $b^2z,  bxz, axz, z^3, b^2x, xz^2, b^2xz, xz^3\}$
\end{tabular}
\caption{\label{subsemis} Subsemigroups ${\mathcal S}(U)$ of ${\mathcal Q}_{{\bf 0.123}}$ for various choices of subsets $U \subseteq H_{12}$.
Boxed elements $s$ additionally satisfy the equation $\omega = xb = s \cdot \Phi(U)$ (see equation (\ref{divispossibility}) in
Theorem 2, and also Figure \ref{sols}).}
\end{center}
\end{figure}

\begin{figure}[h]
\begin{center}
\begin{tabular}{ccc}
 & & Values $s \in {\mathcal Q}_{{\bf 0.123}}$ \\ 
$U$ & $\Phi(U)$ & with $xb = s \cdot \Phi(U)$ \\ \hline

$\{h_4\}$ & $z$  & $\{bxz\}$ \\
$\{h_8\}$ & $a$  & $\{bxz\}$ \\
$\{h_9\}$ & $b$ & $\{x,xz^2,axz\}$ \\
$\{h_{10}\}$ & $x$ & $\{b\}$ \\
$\{h_4,h_8\}$ & $za$ & $\{bx\}$ \\
$\{h_4,h_9\}$ & $zb$ & $\{xz,ax,xz^3\}$ \\
$\{h_4,h_{10}\}$ & $zx$ & $\{bz\}$ \\
$\{h_8,h_9\}$ & $zb$ & $\{xz,ax,xz^3\}$ \\
$\{h_8,h_{10}\}$ & $ax$ & $\{bz\}$ \\
$\{h_9,h_{10}\}$ & $bx$ & $\{\framebox{$1$},z^2,az\}$ \\
$\{h_4,h_8,h_9\}$ & $b$ & $\{x,xz^2,axz\}$ \\
$\{h_4,h_8,h_{10}\}$ & $zax$ & $\{b\}$ \\
$\{h_4,h_9,h_{10}\}$ & $bxz$ & $\{\framebox{$z$},a\}$ \\
$\{h_8,h_9,h_{10}\}$ & $bxz$ & $\{z,\framebox{$a$}\}$ \\
$\{h_4,h_8,h_9,h_{10}\}$ & $bx$ & $\{\framebox{$1$},\framebox{$z^2$},\framebox{$az$}\}$ \\

\end{tabular}
\caption{\label{sols} Values $\Phi(U)$ and all solutions of the equation $xb = s \cdot \Phi(U)$ in
${\mathcal Q}_{{\bf 0.123}}$ for various choices of subsets $U \subseteq H_{12}$.  Boxed elements
are additionally members of ${\mathcal S}(U)$; see Figure \ref{subsemis}.}
\end{center}
\end{figure}

Figure \ref{subsemis} shows the subsemigroups ${\mathcal S}(U)$ for each choice of $U$.  Figure \ref{sols} shows 
all solutions to the equation $\omega = xb = s \cdot \Phi(U)$ in
${\mathcal Q}_{\bf 0.123}$.  Intersecting each set in the
second column of Figure \ref{subsemis} with the corresponding set in third column of Figure \ref{sols}, we obtain the solution
pairs $s$, $U$ for equation (\ref{divispossibility}) in Theorem 2. There are six cases of interest, and they're shown using boxes in 
Figure \ref{subsemis} and Figure \ref{sols}.  
Figure \ref{winningmoves} shows the necessary winning move required by Theorem 2, for each of the six cases.

\begin{figure}[h]
\begin{center}
\begin{tabular}{ccccccc}
$U$ & $s$ & $\Phi(U)$ & $h_i$ & $h_i \rightarrow t$ & $s \cdot \Phi(U)$ & $s \cdot \Phi(t) \cdot \Phi(\frac{\partial}{\partial h_i} P(U))$ \\ \hline
$\{h_9, h_{10}\}$ & \framebox{$1$} &  $bx$ & $h_{10}$ & $h_{10} \rightarrow h_8$ & $bx$ & $1 \cdot a \cdot b = ab = zb$  \\
$\{h_4,h_9,h_{10}\}$ & $\framebox{$z$}$ & $bxz$ & $h_{4}$ &$h_{4} \rightarrow h_1$ & $bx$ & $z \cdot x \cdot bx = zb$ \\
$\{h_8,h_9,h_{10}\}$ & $\framebox{$a$}$ & $bxz$ & $h_8$ &$h_8 \rightarrow h_5$ & $bx$ & $a \cdot x \cdot bx = ab = zb$ \\
$\{h_4,h_8,h_9,h_{10}\}$ & $\framebox{$1$}$ & $bx$ & $h_4$ &$h_4 \rightarrow h_1$ & $bx$ & $1 \cdot x \cdot bxz = zb$ \\
$\{h_4,h_8,h_9,h_{10}\}$ & $\framebox{$z^2$}$ & $bx$ & $h_4$ &$h_4 \rightarrow h_1$ & $bx$ & $z^2 \cdot x \cdot bxz = bz^3 = zb$ \\
$\{h_4,h_8,h_9,h_{10}\}$ & $\framebox{$az$}$ & $bx$ & $h_4$ &$h_4 \rightarrow h_1$ & $bx$ & $az \cdot x \cdot bxz = a^2x^2ab = zb$ \\
\end{tabular}
\caption{\label{winningmoves} Verification of some winning moves for positions of type $\omega = xb$, as
required by Theorem 2.  In each row, the final two columns together
form a move translate $(w,P\mbox{-position})$ (cf equation (\ref{pwinning}) and Figure \ref{NtoPfig}).}
\end{center}
\end{figure}

The verification computations for the other fourteen $N$-position types in ${\bf 0.123}$ are similar.

\section{\normalsize Generalizing Guy \& Smith to misere play}
Having dispensed with correctness to heap size $n$, we're ready to generalize the theorem of
Guy and Smith to misere play.

{\bf Theorem 3}  
Suppose a misere impartial game $\Gamma$ is defined by a finite octal, having $P$
places after the point.  If we can empirically find a misere quotient semigroup ${\mathcal Q}$,
associated outcome partition
${\mathcal Q} = P \bigcup N$,
pretending function $\Phi$, and positive integers
$p$ and $r_0$ such that the equation
\[\Phi(h_{r+p}) = \Phi(h_r)\] holds and the analysis
is correct to heap size $r$ for all $r$ in the range $r_0 \leq r < 2r_0+p+P$,
then it is also correct to heap size $r$ for all $r \geq r_0$, so that the $\Phi$-function has ultimate
period $p$.

{\it Proof.}\/
Because the verification algorithms described in the previous section 
depend only upon move translates and the images of the pretending function on single heap positions,
it suffices to show that the move translates $T_{r+p}$ to heap size $r+p$ are identical to the move translates $T_{r}$,
ie
\[T_{r+p} = T_{r}\]
for all $r \geq 2r_0+p+P$.
For this, the original proof of Guy and Smith for normal play carries over with only minor modifications.
We shamelessly duplicate their language and notation in what follows.

Suppose that the given periodicity relationship has been shown to be correct to heap size $r$ 
for $r_0 \leq r < r_1$,
where $r_1 \geq 2r_0+p+P$.  We want to show it is also true for $r=r_1$.  Let
$h_{r_1+p} \rightarrow h_{s'}h_{t'}$ be a typical move from the heap of size 
$h_{r_1+p}$ involving removing $c$ beans from the heap.  Then
\[(\Phi(h_{r_1+p}),\Phi(h_{s'})\Phi(h_{t'}))\]
is a typical move pair and 
\[(u \Phi(h_{r_1+p}),u \Phi(h_{s'})\Phi(h_{t'}))\]
is a typical move translate with $u \in {\mathcal Q}$ arbitrary.  Then
\[c+s'+t' =  r_1 +p,\]
and $c \leq P$, and we can suppose $s' \leq t'$ so that
\[P + t' + t' \geq r_1+p \geq 2r_0+2p + P,\]
so that
\[t'-p \geq r_0 \geq 0,\]
and therefore $\Phi(h_{t'-p}) = \Phi(h_{t'})$ by the inductive hypothesis.  But there is also a permissible
move $h_{r_1} \rightarrow h_{s'}h_{t'-p}$ so we see that
\[(u \Phi(h_{r_1+p}),u \Phi(h_{s'})\Phi(h_{t'})) = (u \Phi(h_{r_1}),u \Phi(h_{s'})\Phi(h_{t'-p}))\]
are identical move translates.  We've shown that $T_r \subseteq T_{r+p}$.   A similar argument shows
that $T_{r+p} \subseteq T_{r}$.    So $T_r = T_{r+p}$ and the verification algorithms of the previous
section will succeed.

\section{\normalsize Canonical forms \& genera vs misere quotients}
In this section, we contrast the quotient semigroup approach to ${\bf 0.123}$ to two traditional
approaches to misere games:  {\em canonical forms} and {\em genus values}.
\label{canonical}

\subsection{\normalsize Misere canonical forms}

In normal play, the Sprague-Grundy theory describes how to determine the outcome of a sum $G+H$ of two games $G$
and $H$ by computing {\em canonical forms} for each summand---these turn out to be {\em nim-heap 
equivalents}.  We then can imagine that we're playing Nim on each summand $G$ and $H$ instead, and can use {\em nim addition}
to determine the outcome of the sum $G+H$.  At the center of the Sprague-Grundy theory is the equation $G+G = 0$, 
which always holds for an arbitrary normal play combinatorial game $G$.

In misere play, canonical forms can be computed also, but the resulting positions are not nim-heaps in general.
Instead, the canonical form of a typical misere game looks like a complicated tree of options (see \cite{onag}, for example,
in its chapter 12, ``How to Lose when you Must;'' or \cite{ww2}, its chapter 13, ``Survival in the Lost World.'')
Figure \ref{distinguish} illustrates such a tree.
The rules for the reduction of a general misere game to its canonical form game tree are also described in
\cite{onag} and \cite{ww2}.  Nevertheless,
the canonical form 
viewpoint on misere play does not turn out to be so useful in solving wild misere games.  In this section we would like
to shed some light on why this is the case. 

Consider the misdeeds of $\ast 2$, the humble nim heap of size 2.  In normal play, we always have the equation 
\[\ast 2 + \ast 2 = \ast 0 = 0.\]  But in misere play of Nim,
\[\ast 2 + \ast 2 \neq 0.\]  
The left-hand side is a P position, and the right-hand side, an N position.  

It is true in misere Nim that 
\[ (\ast 2 + \ast 2 + \ast 2)  \ \rho \ \ast 2, \]  
ie, these two sums are indistinguishable {\em in misere Nim}.
So we can think of 
\[ \ast 2 + \ast 2 + \ast 2 \mapsto \ast 2, \] 
as a valid outcome-preserving simplication rule in misere Nim.  
The same indistinguishability relation
holds even in a general sum of {\em tame} misere games.  But this doesn't change the fact that
\[\ast 2 + \ast 2 + \ast 2 \neq \ast 2,\] 
since the misere canonical forms of $\ast 2 + \ast 2 + \ast 2$ and $\ast 2$ are {\em not} identical.

It gets worse.  When we move beyond tame games to {\em wild} misere games
such as {\bf 0.123}, we may find that only a weaker indistinguishability relation such as
\begin{equation}
\label{123-eq2}
(\ast 2 + \ast 2 + \ast 2 + \ast 2) \ \rho \ (\ast 2 + \ast 2)
\end{equation}  
is valid.   Such facts are far from obvious or easily proved when one is working 
in the context
of canonical forms.  And for other games---such as ${\bf 4.7}$---there 
is provably no pair of integers $m > n$ such that
\[ \underbrace{\ast 2 + \cdots + \ast 2}_{m \mbox{ copies}} \mapsto \underbrace{\ast 2 + \cdots + \ast 2}_{n \mbox{ copies}} \]
is an outcome-preserving simplification rule \cite{plambeck}.  

So---even the relatively unassuming game $\ast 2$ is an extremely changeable animal 
whose behaviour in sums very much depends on the ``local context'' of the rules of the 
game under
analysis.  More complicated misere games often have even worse behavior, and
involve extremely complicated canonical forms.  
Calculating them explicitly tends to exhaust a computer's memory rapidly. 
By instead restricting our attention to a particular game's 
position sums and concentrating on its quotient rather than its canonical forms, 
it becomes possible to uncover locally valid 
simplication rules such as equation (\ref{123-eq2}), and make 
progress analyzing wild misere games.

\subsection{\normalsize Genera in {\bf 0.123}}
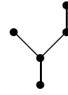
\begin{figure}[h]
\begin{center}
\begin{picture}(45,40)(0,0)

\put(30,0){\circle*{3}}

\put(30,10){\circle*{3}}
\put(30,10){\line(0,-1){10}}

\put(20,20){\circle*{3}}
\put(20,20){\line(1,-1){10}}

\put(40,20){\circle*{3}}
\put(40,20){\line(-1,-1){10}}

\put(40,30){\circle*{3}}
\put(40,30){\line(0,-1){10}}
\end{picture}
\end{center}
\caption{\label{dot123heap6}The heap $h_{6}$ in {\bf 0.123} is $2_+ =\{2\}$, a game of
genus $0^{02}$.}
\end{figure}

\begin{figure}[h]
\begin{center}
\begin{picture}(70,60)(0,0)

\put(35,0){\circle*{3}}

\put(25,10){\circle*{3}}
\put(25,10){\line(1,-1){10}}

\put(20,20){\circle*{3}}
\put(20,20){\line(1,-2){5}}

\put(20,30){\circle*{3}}
\put(20,30){\line(0,-1){10}}

\put(10,40){\circle*{3}}
\put(10,40){\line(1,-1){10}}

\put(30,40){\circle*{3}}
\put(30,40){\line(-1,-1){10}}

\put(30,50){\circle*{3}}
\put(30,50){\line(0,-1){10}}

\put(35,20){\circle*{3}}
\put(35,20){\line(-1,-1){10}}

\put(35,30){\circle*{3}}
\put(35,30){\line(0,-1){10}}

\put(50,30){\circle*{3}}
\put(50,30){\line(0,-1){10}}

\put(40,40){\circle*{3}}
\put(40,40){\line(1,-1){10}}

\put(60,40){\circle*{3}}
\put(60,40){\line(-1,-1){10}}

\put(60,50){\circle*{3}}
\put(60,50){\line(0,-1){10}}

\put(50,20){\circle*{3}}
\put(50,20){\line(0,-1){10}}

\put(50,10){\circle*{3}}
\put(50,10){\line(-3,-2){15}}

\put(65,20){\circle*{3}}
\put(65,20){\line(-3,-2){15}}
\end{picture}
\end{center}
\caption{\label{dot123heap11}The heap $h_{11}$ in {\bf 0.123}, also of
genus $0^{02}$.}
\end{figure}
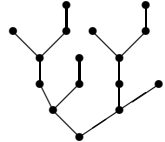

Let the symbol $h_n$ stand for the heap of size $n$ in ${\bf 0.123}$.
The canonical form game trees for $h_6$ and $h_{11}$ are shown 
in Figures \ref{dot123heap6} and \ref{dot123heap11}.    
It's apparent that these trees are not isomorphic---they represent different canonical forms.
If two canonical forms are different,    
there has to be a game $T$ in the global semigroup of all impartial misere combinatorial games
that distinguishes between them \cite{onag}.
But exactly what game $T$ would distinguish between them?
Such a $T$ cannot be a position of ${\bf 0.123}$---if it were, our identification
of these two heaps with the semigroup element $b^2$ (implicit in Figure \ref{0123-semi}) would be invalid.


\begin{figure}[h]
\begin{center}
\begin{picture}(220,55)(-25,0)

\put(100,0){\circle*{3}}

\put(40,10){\circle*{3}}
\put(40,10){\line(6,-1){60}}

\put(20,30){\circle*{3}}
\put(20,30){\line(1,-1){20}}

\put(10,30){\circle*{3}}
\put(10,30){\line(1,0){10}}

\put(10,40){\circle*{3}}
\put(10,40){\line(1,-1){10}}

\put(0,40){\circle*{3}}
\put(0,40){\line(1,0){10}}

\put(30,40){\circle*{3}}
\put(30,40){\line(-1,-1){10}}

\put(30,50){\circle*{3}}
\put(30,50){\line(0,-1){10}}

\put(20,45){\circle*{3}}
\put(20,45){\line(2,-1){10}}

\put(20,55){\circle*{3}}
\put(20,55){\line(2,-1){10}}

\put(50,20){\circle*{3}}
\put(50,20){\line(-1,-1){10}}

\put(50,30){\circle*{3}}
\put(50,30){\line(0,-1){10}}

\put(55,40){\circle*{3}}
\put(55,40){\line(-1,-2){5}}

\put(40,35){\circle*{3}}
\put(40,35){\line(2,-1){10}}

\put(45,45){\circle*{3}}
\put(45,45){\line(2,-1){10}}

\put(80,10){\circle*{3}}
\put(80,10){\line(2,-1){20}}

\put(70,15){\circle*{3}}
\put(70,15){\line(2,-1){10}}

\put(80,20){\circle*{3}}
\put(80,20){\line(0,-1){10}}

\put(85,30){\circle*{3}}
\put(85,30){\line(-1,-2){5}}

\put(70,25){\circle*{3}}
\put(70,25){\line(2,-1){10}}

\put(75,35){\circle*{3}}
\put(75,35){\line(2,-1){10}}

\put(90,25){\circle*{3}}
\put(90,25){\line(-2,-3){10}}

\put(90,35){\circle*{3}}
\put(90,35){\line(0,-1){10}}

\put(80,40){\circle*{3}}
\put(80,40){\line(2,-1){10}}

\put(95,45){\circle*{3}}
\put(95,45){\line(-1,-2){5}}

\put(85,50){\circle*{3}}
\put(85,50){\line(2,-1){10}}

\put(120,10){\circle*{3}}
\put(120,10){\line(-2,-1){20}}

\put(110,10){\circle*{3}}
\put(110,10){\line(1,0){10}}

\put(110,20){\circle*{3}}
\put(110,20){\line(1,-1){10}}

\put(100,20){\circle*{3}}
\put(100,20){\line(1,0){10}}

\put(130,20){\circle*{3}}
\put(130,20){\line(-1,-1){10}}

\put(130,30){\circle*{3}}
\put(130,30){\line(0,-1){10}}

\put(120,25){\circle*{3}}
\put(120,25){\line(2,-1){10}}

\put(120,35){\circle*{3}}
\put(120,35){\line(2,-1){10}}

\put(160,10){\circle*{3}}
\put(160,10){\line(-6,-1){60}}

\put(170,20){\circle*{3}}
\put(170,20){\line(-1,-1){10}}

\end{picture}
\end{center}
\caption{\label{distinguish} The non-{\bf 0.123} position $T = \{2_+3, 2_+20, 3, 1 \}$ 
distinguishes between the ${\bf 0.123}$ positions $h_{6}$ and $h_{11}$ in the global
semigroup of all impartial misere combinatorial games.  But since $T$ never occurs as a position
of ${\bf 0.123}$, the existence of such a $T$ is not relevant to the best play of ${\bf 0.123}$.}
\end{figure}
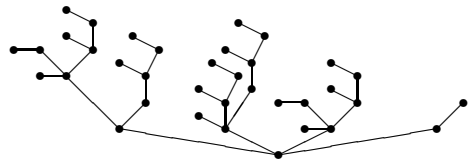

So what does such a game $T$ look like?   One game---obtained via computer search---that 
distinguishes between $h_6$ and $h_{11}$ is
\[ T = \{2_+3, 2_+20, 3, 1 \}, \]
 a game of genus  $0^{20}$.  The game $T$ is illustrated in Figure \ref{distinguish}.
We have
\[\mbox{genus}(h_6+T) = 0^{20}, \mbox{ while} \]
\[\mbox{genus}(h_{11}+T) = 0^{0520}. \]
In particular, the sum $h_6+T$ is a misere N-position, while $h_{11}+T$ is a P-position.
However, the existence of such a $T$ is not relevant to the best play of ${\bf 0.123}$,
for the simple reason that it never occurs in ${\bf 0.123}$!

By way of contrast, consider $h_3+h_3$,
a third position of genus $0^{02}$ in {\bf 0.123}.  Unlike $h_6$ and $h_{11}$,
the position $h_3+h_3$ is of type $z^2$ in ${\bf 0.123}$.  And indeed
the position $h_3+h_3$ {\em can be} distinguished from both $h_6$ and $h_{11}$ in {\bf 0.123}.
For example, $S=h_5+h_9$, a position of type $bx$, distinguishes between $h_6$ and
$h_3 + h_3$:

\[\mbox{genus}(h_6+S) = 0^{02}, \mbox{ while} \]
\[\mbox{genus}(h_{3}+h_{3}+S) = 0^{31}. \]

The former is a misere P-position, while the latter is an N-position.


We conclude that positions that have the same genus symbol may 
or may not be equivalent in the misere quotient semigroup.  

\subsection{\normalsize Summary}
Loosely recapitulating the previous two sections, we might say this:
\begin{quotation}
When we're trying to analyze the misere positions of a game $\Gamma$ with fixed rules,
canonical forms can make positions that are really
the {\em same} in $\mathcal{Q}(\Gamma)$ look {\em different}, while genus computations can
make positions that are really {\em different} in ${\mathcal Q}(\Gamma)$ look the {\em same}.
Only the quotient 
semigroup ${\mathcal Q}(\Gamma)$ precisely 
captures such nuances.
\end{quotation}

\subsection{\normalsize The Sibert-Conway solution to ${\bf 0.77}$ (Misere Kayles)}

The solution to misere Kayles was discovered by William Sibert.  His
original description of it appears in the interesting unpublished manuscript
\cite{sibertorig}.
See also the paper by Sibert and Conway~\cite{sibert}, and the second
edition of Winning Ways~\cite{ww2} (at the end of Volume II).  In this section we
summarize the solution as presented in the paper by Sibert and Conway; in the following one,
we give the corresponding misere quotient.

The PN-positions of Kayles (ie, the positions that a P positions in normal play, but N positions in misere play) are precisely the positions
\[E(5)\ E(4,1),\]
\[E(17,12,9)\ E(20,4,1), \mbox{ or} \] 
\[25 \ E(17,12,9)\ D(20,4,1).\]

The NP-positions (N normal, P misere) are of the form
\[D(5)\ D(4,1),\]
\[E(5)\ D(4,1),\]
\[D(9)\ E(4,1),\]
\[12\ E(4,1), \] 
\[E(17,12,9)\ D(20,4,1), \mbox{ or} \] 
\[25 \ D(9)\ D(4,1).\]

The notation $E(a,b,\ldots)$ (resp. $D(a,b,\ldots)$) refers to any position
composed by taking an even (resp., odd) number of heaps of size $a$ or $b$ or $\ldots$.
For example, 
\[5+4+1+1\]
is a position included in the set 
\[D(5)\ D(4,1)\]
since it is composed by taking a single (ie, odd number) heap of size 5 
and the {\em total number} of 4's or 1's (ie, three) is also odd. 

For every other position in misere Kayles not listed in the
PN- and NP-positions above, its misere and normal play outcomes agree.
\label{kaylessolution}

\begin{figure}[h]
\begin{center}
\begin{tabular}{c|cccccccccccc}
  & 1 & 2 & 3 & 4 & 5 & 6 & 7 & 8 & 9 & 10 & 11 & 12 \\ \hline 
 0+  & $1^{031}$  & $2^{20}$  & $3^{31}$  & $1^{031}$  & $4^{146}$  & \
$3^{31}$  & $2^{20}$  & $1^{13}$  & $4^{046}$  & $2^{20}$  & $6^{46}$  & \
$4^{046}$  \\  
 12+  & $1^{13}$  & $2^{20}$  & $7^{57}$  & $1^{13}$  & $4^{64}$  & $3^{31}$  \
& $2^{20}$  & $1^{031}$  & $4^{64}$  & $6^{46}$  & $7^{57}$  & $4^{64}$  \\  
 24+  & $1^{731}$  & $2^{20}$  & $8^{8[10]}$  & $5^{75}$  & $4^{64}$  & \
$7^{57}$  & $2^{20}$  & $1^{13}$ & \ldots &  &  & 
\end{tabular}
\end{center}
\caption{Genera for {\bf 0.77}}
\end{figure}

\subsection{\normalsize The Kayles misere quotient}
\label{kaylesmisere}

A presentation of the Kayles misere quotient $\mathcal{Q}_{0.77}$ 
can be written down using seven generators

\[\{x, z, w, v, t, f, g \},\]
which first appear in its pretending function at the respective heap sizes 
\[\{1, 2, 5, 9, 12, 25, 27 \}.\]  
Let the symbol $e$ represent the identity of $\mathcal{Q}_{0.77}$. 
The generators satisfy the relations 
\begin{center}
\begin{tabular}{cccc}
$x^2 = e$, & $z^3=z$, & $w^3=w$, & $v^3=v$, \\
\end{tabular}
\end{center}
\begin{center}
\begin{tabular}{cccc}
$t^4 = t^2$, & $f^4=f^2$, & $g^3=g.$ & \\
\end{tabular}
\end{center}

If $\mathcal{Q}_{0.77}$ were isomorphic to the direct product of these seven generators,
it would be a semigroup of order 
$2 \times 3 \times 3 \times 3 \times 4 \times 4 \times 3 = 864$.  
However, many indistinguishability relations
hold between monomials in the generators.   
The actual Kayles quotient is a semigroup of order only 40. Its elements are
shown in Figure \ref{kayles-elements}.   We can
contrast this with the normal play quotient, which is a group of order 16.

\begin{figure}[h]
\begin{center}
\begin{tabular}{ccccc}
$e$ & $x$ & $z$ & $w$ &  $v$ \\
$t$ & $f$ & $g$ & $xz$ & $xw$ \\
$xv$ & $xt$ & $x\!f$ & $xg$ &  $z^2$ \\
$zw$ & $zg$ & $w^2$ & $w\!f$ &  $wg$ \\
$v^2$ & $vt$ & $v\!f$ & $t\!f$ & $xz^2$ \\
$xzw$ & $xzg$ & $xw^2$ & $xw\!f$ &  $xwg$ \\
$xv^2$ & $xvt$ & $xv\!f$ & $xt\!f$ &  $zwg$ \\
$v^2t$ & $vt\!f$ & $xzwg$ & $xv^2t$ &  $xvt\!f$ \\
\end {tabular}
\end{center}
\caption{\label{kayles-elements} The forty elements of the Kayles misere quotient}
\end{figure}

The pretending function for misere Kayles is shown in Figure \ref{077pretend}.  Its
final twelve values repeat indefinitely.
For the sake of comparison, the normal play nim sequence of Kayles is shown in Figure \ref{077seq}.
It also has period twelve.

\begin{figure}
\begin{center}
\begin{tabular}{c|cccccccccccc}
    &  1   &  2   &  3   &  4   &  5  &   6  &   7   &   8 & 9 & 10 & 11 & 12 \\ \hline
0+  & $x$  & $z$  & $xz$  & $x$ &  $w$ & $xz$ & $z$ & $xz^2$ & $v$ & $z$ & $zw$ & $t$ \\
12+ & $xz^2$ & $z$ & $zwx$ & $xz^2$ & $v^2t$ & $xz$ & $z$ & $xvt$ & $wz^2$ & $zw$ & $zwx$ & $wz^2$ \\
24+ & $f$ & $z$ & $g$ & $xwz^2$ & $wz^2$ & $zwx$ & $z$ & $xz^2$ & $g$ & $zw$& $zwx$ & $wz^2$ \\
36+ & $xz^2$ & $z$ &$xz$ & $xz^2$ & $wz^2$ & $zwx$ & $z$ & $xz^2$ & $g$ & $z$ & $zwx$ & $wz^2$ \\
48+ & $xz^2$ & $z$ &$g$ & $xz^2$ & $wz^2$ & $zwx$ & $z$ & $xz^2$ & $wz^2$ & $z$ & $zwx$ & $wz^2$ \\
60+ & $xz^2$ & $z$ &$g$ & $xz^2$ & $wz^2$ & $zwx$ & $z$ & $xz^2$ & $g$ & $zw$& $zwx$ & $wz^2$ \\
72+ & $xz^2$ & $z$ &$g$ & $xz^2$ & $wz^2$ & $zwx$ & $z$ & $xz^2$ & $g$ & $z$& $zwx$ & $wz^2$ \\
84+ & $xz^2$ & $z$ &$g$ & $xz^2$ & $wz^2$ & $zwx$ & $z$ & $xz^2$ & $g$ & $z$& $zwx$ & $wz^2$ \\
96+ & $\ldots$ & & & & & & & & & & & 
\end{tabular}
\caption{\label{077pretend} The pretending function of misere Kayles.}
\end{center}
\end{figure}

\begin{figure}
\begin{center}
\begin{tabular}{c|cccccccccccc}
    &  1   &  2   &  3   &  4   &  5  &   6  &   7   &   8 & 9 & 10 & 11 & 12 \\ \hline
0+  & $1$  & $2$  & $3$  & $1$ &  $4$ & $3$ & $2$ & $1$ & $4$ & $2$ & $6$ & $4$ \\
12+ & $1$ & $2$ & $7$ & $1$ & $4$ & $3$ & $2$ & $1$ & $4$ & $6$ & $7$ & $4$ \\
24+ & $1$ & $2$ & $8$ & $5$ & $4$ & $7$ & $2$ & $1$ & $8$ & $6$& $7$ & $4$ \\
36+ & $1$ & $2$ & $3$ & $1$ & $4$ & $7$ & $2$ & $1$ & $8$ & $2$ & $7$ & $4$ \\
48+ & $1$ & $2$ &$8$ & $1$ & $4$ & $7$ & $2$ & $1$ & $4$ & $2$ & $7$ & $4$ \\
60+ & $1$ & $2$ &$8$ & $1$ & $4$ & $7$ & $2$ & $1$ & $8$ & $6$& $7$ & $4$ \\
72+ & $1$ & $2$ &$8$ & $1$ & $4$ & $7$ & $2$ & $1$ & $8$ & $2$& $7$ & $4$ \\
84+ & $1$ & $2$ &$8$ & $1$ & $4$ & $7$ & $2$ & $1$ & $8$ & $2$& $7$ & $4$ \\
96+ & $\ldots$ & & & & & & & & & & & 
\end{tabular}
\caption{\label{077seq} The nim sequence of normal play Kayles.}
\end{center}
\end{figure}

The multiplication in $\mathcal{Q}_{0.77}$ is given by the Knuth-Bendix
rewriting system shown in Figure \ref{077kb}.  It was 
calculated using the computer algebra package GAP4
\cite{GAP4}.
\begin{figure}
\begin{verbatim}
Rewriting System for Semigroup( [ e, x, z, w, v, t, f, g ] ) with rules 
  [ [ x^2, e ], [ z*v, z*w ], [ z*t, z*w ], [ w*t, z^2 ], 
  [ v*w, z^2 ], [ t^2, v*t ], [ f^2, z^2 ], [ f*g, x*g ], 
  [ v*g, w*g ], [ t*g, w*g ], [ g^2, z^2 ], [ z^3, z ],   
  [ z*w^2, z ], [ w^3, w ],   [ v^3, v ],   [ v^2*f, f ], 
  [ z^2*g, g ], [ w^2*g, g ], [ z*f, x*z ], [ z^2*w, x*w*f ], 
  [ w^2*f, x*z^2 ] ]
\end{verbatim}
\caption{\label{077kb} A Knuth-Bendix rewriting system for misere Kayles ({\bf 0.77}).}
\end{figure}


\begin{figure}
\begin{center}
\begin{picture}(45,140)(0,0)
\setlength{\unitlength}{.8\unitlength}

\put(0,0){$z^2$}
\put(-5,10){\line(-1,1){20}}
\put(19,11){\line(1,1){60}}
\put(-40,40){$vt$}
\put(-35,55){\line(0,1){60}}
\put(-40,120){$v^2$}
\put(-25,135){\line(1,1){20}}
\put(80,80){$w^2$}
\put(75,95){\line(-1,1){60}}
\put(4,160){$e$}

\end{picture}
\caption{\label{077order}The natural partial ordering of the five idempotents of the misere Kayles ({\bf 0.77}) quotient.}
\end{center}
\end{figure}
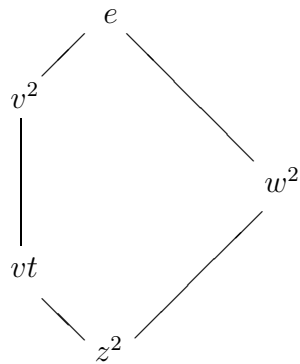


The semigroup $\mathcal{Q}_{0.77}$ has nine pairwise distinguishable P-position types
\[\{x,v,t,z^2,xw,xw^2,xv^2,xvt,xv\!f\}\]
and five idempotents
\[\{e, z^2, w^2,v^2, vt\}.\]

Figure \ref{kaylesmaximalsubgroups} and Figure \ref{077miznorm} together
identify the five maximal subgroups of $\mathcal{Q}_{0.77}$.  There is
a maximal subgroup corresponding to each idempotent (recall section \ref{divisibility}, above).

The natural partial ordering of idempotents is shown in Figure \ref{077order}.

\begin{figure}
\begin{center}
\begin{eqnarray*}
I(e) & = & \{e, x\} \\
I(w^2) & = & \{w^2,w,wx,w^2x\} \\
I(v^2) & = & \{v^2,v,vx,v^2x\} \\
I(vt) & = & \{vt,v^2t,xvt,xv^2t\} \\
I(z^2) & = & \mbox{See Figure \ref{077miznorm}.} \\
\end{eqnarray*}
\caption{\label{kaylesmaximalsubgroups} The five maximal subgroups $I()$ of $\mathcal{Q}_{0.77}$. 
Each subgroup includes precisely those elements that both divide and are divisible by 
an idempotent element that acts as that subgroup's identity.}

\end{center}
\end{figure}

The maximal subgroup $I(z^2)$ is isomorphic to the normal play ${\bf 0.77}$ quotient.  Figure
\ref{077miznorm} shows this correspondence.  

\begin{figure}
\begin{center}
{\small
\begin{tabular}{cccccccc}
$\ast 0$ & $\ast 1$ & $\ast 2$ & $\ast 3$ & $\ast 4$ & $\ast 5$ & $\ast 6$ & $\ast 7$ \\ \hline
$z^2$ & $xz^2$ & $z$ & $xz$ & $xw\!f$ & $w\!f$ & $wz$ & $wxz$ \\
& & & & & & &
\end{tabular} 

\begin{tabular}{cccccccc}
$\ast 8$ & $\ast 9$ & $\ast 10$ & $\ast 11$ & $\ast 12$ & $\ast 13$ & $\ast 14$ & $\ast 15$ \\ \hline
$g$ & $gx$ & $gz$ & $gxz$ & $gw$ & $gwx$ & $gwz$ & $gwxz$
\end{tabular} 
}
\caption{\label{077miznorm}The maximal subgroup $I(z^2)$ of the Kayles misere quotient semigroup 
$\mathcal{Q}_{0.77}$ includes precisely those elements that both divide and are divisible by the
idempotent element $z^2$.  The sixteen elements of 
$I(z^2)$ form a subgroup of the semigroup $\mathcal{Q}_{0.77}$ 
that is isomorphic to the Kayles 
normal play quotient ${ Z}_2 \times { Z}_2 \times { Z}_2 \times { Z}_2$.}
\end{center}
\end{figure}

\section{\normalsize Other games \& discussion}

The techniques described in this paper yield complete analyses for many previously unsolved wild octal games; in
subsequent work we hope to provide a census of such results.   Some examples are available now
at the web site \cite{mg}.  Even where the semigroup techniques fail 
to yield a complete analysis,
the author has observed them to be a powerful tool in extending the analysis of impartial games in misere play.
For each of the many complete analyses of normal play impartial games in the literature, the corresponding
misere quotient calls out to be discovered. 

Commutative semigroup theory is a rich area of mathematical research, and much of it is directly 
applicable in misere play analyses.  
It seems quite likely that
much more can be said in answer to question asked in \cite{ww2} (pg 451):

\begin{quotation}
Are misere analyses really so difficult?
\end{quotation}

\section*{\normalsize Acknowledgement}

I thank Aaron Siegel and Dan Hoey for drawing my attention to errors in section \ref{kaylesmisere} of the first draft
of this manuscript.

\end{document}